\documentclass[11pt,a4paper]{article}
\usepackage{amsmath}
\usepackage{amsfonts}
\usepackage{amssymb}
\usepackage{graphicx}

\usepackage[english]{babel}

\newcommand{\norm}[1]{\left\Vert #1\right\Vert}

\newcommand{\R}{   {\ifmmode{{\mathbb R}}\else{$\mathbb R$}\fi}}
\newcommand{\N}{   {\ifmmode{{\mathbb N}}\else{$\mathbb N$}\fi}}

\newcommand{\K}{\mathcal{K}}
\newcommand{\demo}{\noindent\textit{Proof. }}

\newcommand{\dt}{\frac{d}{dt}} 
\newcommand{\dto}{\frac{d}{dt}_{|t=0}}

\begin{document}

\title{\bf Stability of uniformly bounded switched systems\\ {\large and }\\ Observability}
\author{Moussa BALDE\footnote{LGA-LMDAN Dept. de Math\'ematiques et Informatique,UCAD, Dakar-Fann, Senegal, E-mail: moussa.balde.math@ucad.edu.sn}, Philippe JOUAN\footnote{LMRS, CNRS UMR 6085, Universit\'e
    de Rouen, avenue de l'universit\'e BP 12, 76801
    Saint-Etienne-du-Rouvray France. E-mail: Philippe.Jouan@univ-rouen.fr}, Sa\"id NACIRI }

\date{\today}

\maketitle

\begin{abstract}
This paper mainly deals with switched linear systems defined by a pair of Hurwitz matrices that share a common but not strict quadratic Lyapunov function. Its aim is to give sufficient conditions for such a system to be GUAS.

We show that this property of being GUAS is equivalent to the uniform observability on $[0,+\infty)$ of a bilinear system defined on a subspace whose dimension is in most cases much smaller than the dimension of the switched system.

Some sufficient conditions of uniform asymptotic stability are then deduced from the equivalence theorem, and illustrated by examples.

The results are partially extended to nonlinear analytic systems.

\vskip 0.2cm

Keywords: Switched systems; Asymptotic stability; Quadratic Lyapunov functions; Observability; Bilinear systems.

\vskip 0.2cm

AMS Subject Classification: 93B07, 93C99, 93D20.

\end{abstract}

\section{Introduction.}
One of the most important issues in the field of switched systems is that of uniform asymptotic stability (see \cite{Liberzon}), and the Lyapunov functions theory is known to be a very powerful tool in this connection. Indeed the existence of a common strict Lyapunov function ensures that the system is GUAS (Globally Uniformly Asymptotically Stable).
Moreover it has been proved in \cite{dw} that a GUAS linear switched system always has a common strict Lyapunov function (this result was extended to nonlinear switched systems in \cite{MAG00}). In \cite{MBC08} Mason, Boscain, and Chitour showed that such a common strict Lyapunov function  can always be chosen to be polynomial. Despite this fact a quadratic one does not always exist (\cite{dw}, \cite{MBC08}).

Among other approaches of the GUAS property, let us quote the worst trajectory method (see  \cite{bbm} and references therein), its characterization via optimal control techniques (\cite{ML06}), and the commutations relations method (see for instance \cite{agr}).

In the present paper we deal first, and mainly, with pairs of matrices that share a common but not strict in general quadratic Lyapunov function, in other words we are interested in the asymptotic properties of systems whose trajectories are known to be bounded (by a quadratic function in the linear case). Properties of switched systems in this framework have already been studied in \cite{RSD10}, \cite{SVR}, \cite{BJ11} (similar properties of nonlinear sytems are stated in \cite{BM05} and \cite{JN13}). In particular our paper follows \cite{BJ11} and makes one of its results much more precise.

To this end we propose a new method that consists in showing that the property of being GUAS is equivalent to the uniform observability on $[0,+\infty)$ of a bilinear system defined in a natural way on a subspace whose dimension is in most cases much smaller than the dimension of the state space of the switched system.

Sufficient conditions of uniform asymptotic stability are then deduced from this observability property. To our knowledge, these conditions are new. They are moreover checkable as shown by the examples of Section \ref{exemples}.

Our results are stated only for pairs of matrices, but can be rather easily extended to finite families of matrices that share a weak quadratic Lyapunov function. However the statement of the result is in that case quite complicated (see Conclusion for more details) and does not introduce new ideas. For these reasons our discussion is restricted to pairs of matrices.

Some relationships between observability and asymptotic stability of linear or nonlinear switched systems have already been established (see for instance \cite{HLAS05} where the output depends on a family of weak Lyapunov functions, and for the linear case  \cite {H04} which makes use of the fact that a symmetric negative matrix $Q$ can be written as $Q=-C^TC$, with $\ker(Q)=\ker(C)$), but to our knowledge, they did not lead to a reduction of the dimension of the state space, which is the main interest of our method.
\vskip 0.2cm

Let us now state our results. We deal with a pair $\{B_0,B_1\}$ of $d\times d$ Hurwitz matrices that share a common, but not strict in general, quadratic Lyapunov function.

We can assume without loss of generality that the Lyapunov matrix is the identity, so that the matrices $B_i$ verify $B_i^T+B_i\leq 0$ for $i=0,1$.

In this setting the linear subspace
$$
\K=\K_0\bigcap \K_1, \quad\mbox{where}\quad \K_i=\ker(B_i^T+B_i)\quad i=0,1,
$$
turns out to be a fundamental object (see \cite{BJ11}, \cite{SVR}, and \cite{JN13} for the nonlinear case).

In the previous paper \cite{BJ11} it was proved that the system is GUAS as soon as $\K=\{0\}$. However this condition is not necessary, and it is possible to build GUAS systems for which $\dim \K=d-1$ regardless of $d$ (see Example \ref{dmoinsun}).

Firstly we show that a necessary and sufficient condition for the switched system to be GUAS is that a certain bilinear system in $\K$ is observable on $[0,+\infty)$. More accurately the convexified matrix $B_{\lambda}=(1-\lambda)B_0+\lambda B_1$ is given, according to the decomposition $\R^d=\K\oplus \K^{\bot}$, by:
$$
B_{\lambda}=\begin{pmatrix}A_{\lambda}&-C_{\lambda}^T\\C_{\lambda}&D_{\lambda}\end{pmatrix},
$$
with $A_{\lambda}^T+A_{\lambda}=0$ for $\lambda\in [0,1]$, and $D_{\lambda}^T+D_{\lambda}<0$ for $\lambda\in (0,1)$.

Consider the bilinear system in $\K$ with output in $\K^{\bot}$:
$$
(\Sigma)=\ \left\{
\begin{array}{rl}
\dot{x} & =A_{\lambda} x\\
y & =C_{\lambda} x
\end{array}
\right.
\quad \lambda \in [0,1],\quad x\in \K,\quad  y\in \K^{\bot}.
$$
Our main result is then:

\vskip 0.2cm

\noindent{\bf Theorem \ref{CNS}}

\textit{The switched system is GUAS if and only if $(\Sigma)$ is observable for all inputs $\lambda$: $[0,+\infty)\mapsto[0,1]$.}

\vskip 0.2cm

The interest of this equivalence is twofold: on the one hand the bilinear observed system is in $\K$, whose dimension is most often (but not always) much smaller than $d$. On the other hand the observability of $(\Sigma)$ is rather easier to deal with than the asymptotic properties of the switched system.

Thanks to Theorem \ref{CNS} we obtain the following sufficient conditions of uniform asymptotic stability (the definition of the set $G$ is rather technical and the reader is referred to Section \ref{Observabilite}):

\vskip 0.2cm

\noindent{\bf Theorem \ref{Final}}

\textit{The switched system is GUAS as soon as the pair $(C_{\lambda},A_{\lambda})$ is observable for every $\lambda \in [0,1]$, and one of the following conditions holds:
\begin{enumerate}
	\item the set $G$ is discrete;
	\item $\dim \K\leq 2$.
\end{enumerate}
In particular the switched system is GUAS if $\ker C_{\lambda}=\{0\}$ for $\lambda \in[0,1]$.}

\vskip 0.2cm

As we know of no system which is not GUAS and such that the pair $(C_{\lambda},A_{\lambda})$ is observable for every $\lambda \in [0,1]$ we make the following conjecture:

\vskip 0.2cm

\noindent{\bf Conjecture}

\textit{The switched system is GUAS if and only if the pair $(C_{\lambda},A_{\lambda})$ is observable for every $\lambda \in [0,1]$.}

\vskip 0.1cm

This conjecture, which is equivalent to saying that the matrix $B_{\lambda}$ is Hurwitz for every $\lambda \in [0,1]$, is discussed in Section \ref{concluding}.

\vskip 0.2cm

From \cite{dw} we known that a GUAS switched system has always a common strict Lyapunov function.
One might think that under the assumption of the existence of a common non strict quadratic Lyapunov function, a GUAS system has a common strict quadratic Lyapunov function. This conjecture is wrong, as shown by Example \ref{Paolo}, due to Paolo Mason.

\vskip 0.2cm

The same method is then applied to nonlinear switched systems. For several reasons the linear statements are only partially extended to the nonlinear case, though the systems are required to be analytic. On the first hand the set $\K$ is no longer a vector subspace, not even a submanifold, neither connected, in general. On the second one no inner product is naturally related to general Lyapunov functions, so that it is not possible to get a decomposition of the vector fields as nice as the previous one. To finish the bilinear observed system we deal with in the linear case is observable if and only if the origin is distinguishable from the other points. Of course this property does not hold for general observed systems.

\vskip 0.2cm

The paper is organized as follows: Section \ref{linear} is devoted to the main linear result, the equivalence between asymptotic stability and observability of the bilinear observed system. In Section \ref{observe} the behaviour of the observed system is analyzed, stability statements are deduced, and linear examples are provided. Section \ref{nonlinear} deals with the nonlinear case.


\section{Asymptotic stability of linear systems}\label{linear}

\subsection{Preliminaries}\label{Preliminaries}

We will use the symbol $X$ to denote the elements of $\R^d$. When $\R^d$ will be decomposed as $\R^d=\K\oplus \K^{\bot}$ we will write
$$
X=x+y \ \mbox{ or }\  X=\begin{pmatrix}x\\y\end{pmatrix} \quad\mbox{where}\quad x\in \K, \ \ y\in \K^{\bot}.
$$
As explained in the introduction we deal with a pair $\{B_0,B_1\}$ of $d\times d$ Hurwitz matrices, assumed to share a common, but not strict in general, quadratic Lyapunov function. More accurately there exists a symmetric positive definite matrix $P$ such that the symmetric matrices $B_i^TP+PB_i$ are nonpositive ($B^T$ stands for the transpose of $B$).
Since the Lyapunov matrix $P$ is common to the $B_i$'s we can assume without loss of generality that $P$ is the identity matrix, in other words that $B_i^T+B_i$ is nonpositive definite for $i=0,1\ $:
\begin{equation}\label{condition}
\forall X\in \R^d, \qquad X^T(B_i^T+B_i)X\leq 0 \qquad \mbox{for}\ \ i=0,1.
\end{equation}

\noindent{\bf Norms}. The natural scalar product of $\R^d$ in this context is the canonical one, defined by $\left\langle X,Y\right\rangle=X^TY$ (it would be $X^TPY$ if the Lyapunov matrix were $P$). The norm of $\R^d$ is consequently chosen to be $\norm{X}=\sqrt{X^TX}$.

\vskip 0.2cm

\noindent{\bf The switched system}. We consider the switched system in $\R^d$
\begin{equation}\label{ssystem}
\dot{X}=\dt X=B_{u(t)}X
\end{equation}
where the input, or switching law, $t\longmapsto u(t)$ is a measurable function from $[0,+\infty)$ into the discrete set $\{0,1\}$.

Such a switching law being given, the solution of (\ref{ssystem}) for the initial condition $X$ is known to be
$$
t\longmapsto \Phi_u(t)X,
$$
where $t\longmapsto \Phi_u(t)$ is the solution of the matrix equation $\dot{M}=B_{u(t)}M$, $\Phi_u(0)=I_d$, or in integral form:
\begin{equation}\label{Int}
\Phi_u(t)=I_d+\int_0^tB_{u(s)}\Phi_u(s)\ ds.
\end{equation}

\noindent{\bf The $\omega$-limit sets}.
For $X\in \R^d$ we denote by $\Omega_u(X)$ the set of $\omega$-limit points of $\{\Phi_u(t)X;\ t\geq 0\}$, that is the set of limits of sequences $(\Phi_u(t_j)X)_{j\geq 0}$, where $(t_j)_{j\geq 0}$ is strictly increasing to $+\infty$.

Thanks to Condition (\ref{condition}), the norm $\norm{\Phi_u(t)X}$ is nonincreasing, and $\Omega_u(X)$ is a compact and connected subset of a sphere $\mathcal{S}(r)=\{x\in \R^d;\ \norm{x}=r\}$ for some $r\geq 0$ (see Proposition 1 of \cite{BJ11}).

\newtheorem{GUAS}{Definition}
\begin{GUAS}\label{GUAS}
The switched system is said to be Globally Uniformly Asymptotically Stable, or GUAS in short, if for every switching law $u$ the system (\ref{ssystem}) is globally asymptotically stable, that is
$$
\forall X\in\R^d \qquad \Phi_u(t)X\longrightarrow_{t\rightarrow +\infty}0.
$$
\end{GUAS}


\subsection{Hurwitz Property and Observability of linear systems}\label{Hurwitzmatrix}

\newtheorem{cnsHurwitz}{Theorem}
\begin{cnsHurwitz} \label{cnsHurwitz}
Let $B$ be a $d\times d$-matrix that satisfies $B^T+B\leq 0$ and let $\K=\ker(B^T+B)$. Up to an orthogonal transformation and according to the orthogonal decomposition $\R^d=\K\oplus \K^{\bot}$, one has
$$
B=\begin{pmatrix}A&-C^T\\C&D\end{pmatrix},
$$
with $A^T+A=0$ and $D^T+D<0$.

Then $B$ is Hurwitz if and only if the pair $(C,A)$ is observable.
\end{cnsHurwitz}

\demo

This theorem can be proved in a classical way by applying Theorem 3.8 of \cite{TSH01} to the pair $(B,B^T+B)$. However our proof is closer to the one of the forthcoming main result, and in our opinion enlights this last.

Let $(b_1,\dots,b_d)$ be an orthonormal basis of $\R^d$ such that $(b_1,\dots,b_k)$ span $\K$. In that basis
$$
B=\begin{pmatrix}A&C_1\\C&D\end{pmatrix}
$$
according to the decomposition $\R^d=\K\oplus \K^{\bot}$. The condition $B^T+B\leq 0$ being invariant under orthogonal transformations, and  $\K$ being equal to $\ker(B^T+B)$, we obtain at once $A^T+A=0$, $C_1=-C^T$ and $D^T+D<0$.

Let us now consider the observed linear system in $\K$ with output in $\K^{\bot}$:
$$
(\Sigma)= \quad \left\{
\begin{array}{ll}
\dot{x} & =Ax\\
y & =Cx
\end{array}
\right.
$$

If $(\Sigma)$ is not observable, then there exists $x\in \K$, $x\neq 0$, such that $Ce^{tA}x=0$ for all $t\in \R$. Therefore
$$
\dt\begin{pmatrix}e^{tA}x\\0\end{pmatrix}=B\begin{pmatrix}e^{tA}x\\0\end{pmatrix}, \quad\mbox{so that}\quad  e^{tB}\begin{pmatrix}x\\0\end{pmatrix}=\begin{pmatrix}e^{tA}x\\0\end{pmatrix}.
$$
The matrix $e^{tA}$ being a rotation one, this vector does not tend to $0$ as $t$ goes to $+\infty$. This shows that $B$ is not Hurwitz.

Conversely let us assume that $B$ is not Hurwitz: the real parts of its eigenvalues being nonpositive there exists a periodic, or constant, non zero trajectory  $e^{tB}X$. Since the function $t\longmapsto\norm{e^{tB}X}$ is non increasing, it is constant, so that:
$$
\forall t\in\R, \ \ 0=\dt \norm{e^{tB}X}^2=(e^{tB}X)^T(B^T+B)e^{tB}X,
$$
which shows that the whole trajectory $e^{tB}X$, $t\in\R$, is contained in $\K$.

Finally let us write $X=(x, 0)^T$ in $\R^d=\K\oplus \K^{\bot}$. Then
$$
e^{tB}X=\begin{pmatrix}e^{tA}x\\0\end{pmatrix} \quad\mbox{ with } \quad Ce^{tA}x=0,
$$
in other words the output of $(\Sigma)$ does not distinguish between the initial states $x$ and $0$.

\hfill $\Box$


\subsection{Main result}\label{mainresult}

Let $B_0$ and $B_1$ be two $d\times d$-matrices that satisfy $B_i^T+B_i\leq 0$, and let $\K_i=\ker(B_i^T+B_i)$ for $i=0,1$. Notice that $\K_i$ is also the set of $X\in\R^d$ such that $\displaystyle\dto \norm{e^{tB_i}X}^2=0$. Here and subsequently $\K$ stands for
$$
\K=\K_0\bigcap \K_1=\ker(B_0^T+B_0) \bigcap\ker(B_1^T+B_1).
$$

\vskip 0.2cm

For $\lambda\in [0,1]$ we write $B_{\lambda}=(1-\lambda)B_0+\lambda B_1$. Let us firstly state the following easy but useful lemma.

\newtheorem{noyau}{Lemma}
\begin{noyau} \label{noyau}
For all $\lambda \in (0,1)$,  $\ker (B_{\lambda}^T+B_{\lambda})=\K$.

Moreover, up to an orthogonal transformation and according to the orthogonal decomposition $\R^d=\K\oplus \K^{\bot}$, $B_{\lambda}$ is given by
$$
B_{\lambda}=\begin{pmatrix}A_{\lambda}&-C_{\lambda}^T\\C_{\lambda}&D_{\lambda}\end{pmatrix},
$$
with $A_{\lambda}^T+A_{\lambda}=0$ for $\lambda\in [0,1]$, and $D_{\lambda}^T+D_{\lambda}<0$ for $\lambda\in (0,1)$.
\end{noyau}

\demo

We have only to prove that $\ker (B_{\lambda}^T+B_{\lambda})=\K$ for $\lambda\in (0,1)$, the proof of the second assertion being similar to the beginning of the one of Theorem \ref{cnsHurwitz}.
If $X^T(B_{\lambda}^T+B_{\lambda})X=0$ for some $X\in\R^d$, then
$$
0=(1-\lambda)X^T(B_0^T+B_0)X+\lambda X^T(B_1^T+B_1)X.
$$
But $X^T(B_i^T+B_i)X\leq 0$ for $i=0,1$, and since $\lambda\neq 0,1$ we obtain $X^T(B_0^T+B_0)X=X^T(B_1^T+B_1)X=0$, that is $X\in \K_0\cap \K_1$.

The converse is straightforward.

\hfill $\Box$

\vskip 0.2cm

\noindent{\bf Remark}. The strict inequality $D_{\lambda}^T+D_{\lambda}<0$ does not hold for $\lambda=0,1$ whenever $\K_0$ or $\K_1$ is strictly larger than $\K$. However the non strict inequality $D_{\lambda}^T+D_{\lambda}\leq 0$ holds for $\lambda\in [0,1]$.

\vskip 0.5cm

In the same way as in Section \ref{Hurwitzmatrix} we will consider the bilinear controlled and observed system:
$$
(\Sigma)=\ \left\{
\begin{array}{rl}
\dot{x} & =A_{\lambda} x\\
y & =C_{\lambda} x
\end{array}
\right.
$$
where $\lambda \in [0,1]$, $x\in \K$, and $y\in \K^{\bot}$.

\newtheorem{inftyObserve}[GUAS]{Definition}
\begin{inftyObserve}\label{inftyObserve}
The system $(\Sigma)$ is said to be uniformly observable on $[0,+\infty)$ if for any measurable input $t\longmapsto \lambda(t)$ from $[0,+\infty)$ into $[0,1]$, the output distinguishes any two different initial states, that is
$$
\forall x_1\neq x_2\in \K \quad m\{t\geq 0;\ C_{\lambda(t)} x_1(t)\neq C_{\lambda(t)} x_2(t)\}>0,
$$
where $m$ stands for the Lebesgue measure on $\R$, and $x_i(t)$ for the solution of $\dot{x}=A_{\lambda(t)} x$ starting from $x_i$, for $i=1,2$.
\end{inftyObserve}

\noindent{\bf Remarks}
\begin{enumerate}
	\item As the ouput depends explicitly on the input, it is measurable but not necessarily continuous. It is the reason why our definition of observability involves the Lebesgue measure.
	\item The observability on $[0,+\infty)$ is not equivalent to the observability on bounded time intervals (See Section \ref{obsinfty} and Examples \ref{Kdeux}-\ref{dmoinsun}).
	\item The system being linear with respect to the state, it is clearly observable for a given input if and only if the output does not vanish for almost every $t\in[0,+\infty)$ as soon as the initial state is different from $0$.
\end{enumerate}

We are now in a position to state our main result:

\newtheorem{CNS}[cnsHurwitz]{Theorem}
\begin{CNS} \label{CNS}
The switched system is GUAS if and only if $(\Sigma)$ is uniformly observable on $[0,+\infty)$.
\end{CNS}

\demo

Let us first assume that the switched system is not GUAS. There exist a measurable input $t\longmapsto u(t)$ from $[0,+\infty)$ into $\{0,1\}$ and an initial state $X\in\R^d$ for which the switched system does not converge to $0$.

Let $l$, with $\norm{l}=r>0$, be a limit point for $X$, and $(t_j)_{j\geq 0}$ a strictly increasing sequence such that
$$
l=\lim_{j\rightarrow +\infty}\Phi_u(t_j)X.
$$
Let $\tau$ be an arbitrary positive number and let us define the sequence $(\phi_j)_{j\geq 0}$ by $\phi_j(t)=\Phi_u(t_j+t)X$ for $t\in [0,\tau]$.
Each function $\phi_j$ verifies
$$
\phi_j(t)=\phi_j(0)+\int_0^t B_{u(t_j+s)}\phi_j(s)\ ds.
$$
On the other hand the sequence $(B^j)_{j\geq 0}$ of functions from $[0,\tau]$ to $\mathcal{M}(d;\R)$, defined by $B^j(s)=B_{u(t_j+s)}$, is bounded, and consequently converges weakly ($*$-weakly to be precise) in $L^{\infty}([0,\tau],\mathcal{M}(d;\R))$, up to a subsequence that we continue to denote by $(B^j)_{j\geq 0}$.

Moreover the limit takes its values in the convexification of $\{B_0,B_1\}$ (see \cite{Sontag98}, Lemma 10.1.3, page 424), and can be written
$$
B_{\lambda(t)}=(1-\lambda(t))B_0+\lambda(t)B_1,
$$
where $t\longmapsto \lambda(t)$ is a measurable function from $[0,\tau]$ into $[0,1]$.

Let us denote by $\psi$ the absolutely continuous and $\R^d$-valued function defined on $[0,\tau]$ by the equation
$$
\psi(t)=l+\int_0^t B_{\lambda(s)}\psi(s)\ ds.
$$

According to Theorem 1, page 57, of \cite{Sontag98}, the sequence $(\phi_j)_{j\geq 1}$ converges uniformly on $[0,\tau]$ to $\psi$.
Moreover this function takes its values in $\Omega_u(X)$, so that
$$
\forall t\in [0,\tau] \qquad \norm{\psi(t)}^2=\norm{\psi(0)}^2=\norm{l}^2=r^2>0.
$$
Thus we have for almost every $t\in [0,\tau]$
\begin{equation}\label{egalzero}
\dt \norm{\psi(t)}^2=\psi(t)^T(B_{\lambda(t)}^T+B_{\lambda(t)})\psi(t)=0.
\end{equation}

\newtheorem{Inter}[noyau]{Lemma}
\begin{Inter} \label{Inter}
For all $t\in [0,\tau]$ the vector $\psi(t)$ belongs to $\K=\K_0\cap \K_1$. In particular $l=\psi(0)\in \K$.
\end{Inter}

\demo
For almost every $t\in [0,\tau]$ we have
$$
\begin{array}{l}
\psi(t)^T(B_{\lambda(t)}^T+B_{\lambda(t)})\psi(t)\\
\qquad \qquad \qquad  =(1-\lambda(t))\psi(t)^T(B_0^T+B_0)\psi(t)+\lambda(t)\psi(t)^T(B_1^T+B_1)\psi(t).
\end{array}
$$
But according to (\ref{egalzero}) and
$$
\psi(t)^T(B_i^T+B_i)\psi(t)\leq 0 \quad\mbox{ for }\quad i=0,1
$$
we obtain for almost every $t\in [0,\tau]$
$$
\begin{array}{rl}
\lambda(t)\neq 0 & \Longrightarrow \psi(t) \in \K_1 \\
\lambda(t)\neq 1 & \Longrightarrow \psi(t) \in \K_0
\end{array}
$$
so that $\psi(t)\in \K_0\cup \K_1$. Assume that $l\in \K_0\setminus \K_1$. Then for some $T$, $0<T\leq\tau$, we have
$$
\psi([0,T])\cap \K_1 =\emptyset, \quad \mbox{ hence } \quad \psi(t)^T(B_1^T+B_1)\psi(t)<0
$$
for all $t\in [0,T]$. Consequently $\lambda(t)=0$ for almost every $t\in [0,T]$, and
$$
\psi(t)=e^{tB_0}\psi(0).
$$
But according to the Hurwitz property of $B_0$, the norm $\norm{\psi(t)}$ would be strictly decreasing, in contradiction with its belonging to $\Omega_u(X)$. Consequently $l\in \K$, and in the same way $\psi(t)\in \K$ for $t\in [0,\tau]$.

\hfill $\Box$

\vskip 0.2cm
\noindent{\it End of the proof of Theorem \ref{CNS}}.

As $\psi(t)$ is in $\K$ for all $t$ it can be written according to the decomposition $\R^d=\K\oplus \K^{\bot}$ as a column
$$
\psi(t)=\begin{pmatrix}\phi(t)\\ 0 \end{pmatrix}.
$$
Moreover the derivative of $\psi(t)$ is also in $\K$ for almost every $t$. This derivative is
$$
\dt \psi(t)=B_{\lambda(t)}\psi(t)=\begin{pmatrix}A_{\lambda(t)}\phi(t)\\ C_{\lambda(t)}\phi(t) \end{pmatrix}
$$
and the belonging of $\displaystyle \dt \psi(t)$ to $\K$ turns out to be
$$
C_{\lambda(t)}\phi(t)=0 \quad \mbox{ for almost every }\ t\in [0,\tau].
$$
The conclusion is that $\phi$ is a trajectory of
$$
(\Sigma)=\ \left\{
\begin{array}{rl}
\dot{x} & =A_{\lambda} x\\
y & =C_{\lambda} x
\end{array}
\right.
$$
for which the output vanishes almost surely. Notice that $\phi$ does not vanish since $\norm{\phi(t)}^2=\norm{\phi(0)}^2=\norm{l}^2=r^2>0$ for all $t\in [0,\tau]$.

To conclude this part of the proof it remains to notice that $\phi$ can be extended to $[0,+\infty)$ with the same properties: starting from the final point $\psi(\tau)$ we can obtain a similar limit trajectory on $[\tau,\tau_1]$ for any $\tau_1>\tau$.

This proves that $(\Sigma)$ is not uniformly observable on $[0,+\infty)$.

\vskip 0.2cm

Conversely assume the switched system to be GUAS. It is a well known fact that the convexified system is also GUAS (see \cite{MBC08}). If there exists for $(\Sigma)$ an input defined on $[0,+\infty)$ and with values in $[0,1]$ such that the trajectory $\phi(t)$ for the initial condition $\phi(0)\neq 0$ satisfies $C_{\lambda(t)} \phi(t)=0$ for $t\geq 0$ then
$$
\psi(t)=\begin{pmatrix}\phi(t)\\ 0 \end{pmatrix}
$$
is, for the same input, a trajectory of the convexified switched system that does not converge to $0$, a contradiction.

\hfill $\Box$


\section{Applications of observability}\label{observe}

\subsection{Observability of the bilinear system}\label{Observabilite}

We consider now the controlled and observed bilinear system
$$
(\Sigma)=\ \left\{
\begin{array}{rl}
\dot{x} & =A_{\lambda} x\\
y & =C_{\lambda} x
\end{array}
\right.
$$
where $x\in \K$, $y\in \K^{\bot}$ and $\lambda \in [0,1]$. Notice that the matrices $A_0$ and $A_1$ are skew-symmetric, so that the trajectories of $(\Sigma)$ are contained in spheres. We will denote by $\mathcal{S}^{k-1}$, where $k=\dim \K$, the unit sphere of $\K$.

A solution $t\longmapsto x(t)$ of $(\Sigma)$ on $I=[0,T]$ or $I=[0,+\infty)$,  which is in $\mathcal{S}^{k-1}$ and satisfies
$$
C_{\lambda(t)} x(t)=0 \qquad\mbox{ for almost every } t\in I,
$$
will be called a \textbf{bad trajectory} on $I$.

The purpose is to find conditions for $(\Sigma)$ to be uniformly observable on $[0,+\infty)$.
An obvious necessary condition is that $(\Sigma)$ is observable for every constant input, that is the pair $(C_{\lambda},A_{\lambda})$ is observable for every $\lambda \in [0,1]$. Since $B_0$ and $B_1$ are Hurwitz this property is guaranteed for $\lambda =0,1$ (although $\K$ is not necessarily the kernel of $B_i^T+B_i$, $i=0,1$, this can be easily shown using the same kind of arguments as in the proof of Theorem \ref{cnsHurwitz}).

A sufficient condition is that $(\Sigma)$ is uniformly observable on every bounded interval $[0,T]$, $T>0$, that is $(\Sigma)$ is uniformly observable in the usual meaning.

\vskip 0.1cm

\noindent{\bf Remark (*)}. \textit{Under the condition that the pair $(C_{\lambda},A_{\lambda})$ is observable for every $\lambda \in [0,1]$, no bad trajectory can be constant. This remark is used several times in the forthcoming proofs.}

\vskip 0.1cm

The first task is to characterize the locus where the ouput vanishes (Section \ref{badlocus}). Then we will state some sufficient conditions of uniform observability (Section \ref{Uniformobs}), and of uniform observability on $[0,+\infty)$ (Section \ref{obsinfty}).


\subsection{The bad locus}\label{badlocus}

The condition
$$
\exists\ \lambda\in [0,1]\quad \mbox{ such that }\ \ C_{\lambda} x=(1-\lambda) C_0x+\lambda C_1x=0
$$
is equivalent to saying that $C_0x$ and $C_1x$ are colinear and in opposite directions, that last condition being due to $\lambda \in [0,1]$. Consequently the set of points $x\in \K$ for which there exists $\lambda \in [0,1]$ such that $C_{\lambda} x=0$ can be characterized in the following way:
\begin{equation}\label{carac}
\begin{array}{ll}
\exists \lambda \in [0,1] \mbox{ s.t. } C_{\lambda} x=0 & \Longleftrightarrow \left\langle C_0x,C_1x\right\rangle +\norm{C_0x}\norm{C_1x}=0
\end{array}
\end{equation}
Here the scalar product and the norm are the restrictions to $\K^{\bot}$ of the ones of $\R^d$.

Another helpful characterization of this set makes use of the exterior product $C_0x\wedge C_1x$, which can be considered as the $\displaystyle \frac{k'(k'-1)}{2}$ vector of all the $2\times 2$ minors of the $k'\times 2$ matrix $(C_0x\ C_1x)$ (here $k'$ stands for the dimension of $\K^{\bot}$). We obtain
\begin{equation}\label{ExtWedge}
\begin{array}{ll}
\exists \lambda \in [0,1] \mbox{ s.t. } C_{\lambda} x=0 & \Longleftrightarrow C_0x\wedge C_1x=0 \mbox{ and } \left\langle C_0x,C_1x\right\rangle\leq 0.
\end{array}
\end{equation}
The interest of this last characterization lies in the fact that the exterior product being bilinear, the condition $C_0x\wedge C_1x=0$ can be differentiated.

The set of points that satisfy (\ref{carac}) or (\ref{ExtWedge}) will be denoted by $F$. Since Formula (\ref{ExtWedge}) is clearly invariant under multiplication by a constant, this set is a cone.

We also write $N=\ker C_0\bigcap\ker C_1$ and $F_0=F\setminus N$. For $x\in F_0$ it is clear that $C_0x\neq C_1x$ so that the unique $\lambda$ such that $C_{\lambda} x=0$ is given by:
\begin{equation}\label{ldex}
\lambda(x)=\frac{\left\langle C_0x-C_1x,C_0x\right\rangle}{\norm{C_0x-C_1x}^2},
\end{equation}
that is $x\longmapsto\lambda(x)$ is the restriction to $F_0$ of an analytic function defined on $\K\setminus \{C_0x=C_1x\}$.

Any bad trajectory lies in the intersection of $F$ with $\mathcal{S}^{k-1}$, and as long as it does not meet $N$, that is as long as it remains in $F_0$, Formula (\ref{ldex}) shows that $t\longmapsto\lambda(t)$ and $t\longmapsto x(t)$ are analytic (more accurately $t\longmapsto\lambda(t)$ is almost everywhere equal to an analytic function).


\subsection{Uniform observability}\label{Uniformobs}

Let $(\lambda(t),x(t))$ be a bad trajectory on $[0,T]$ for some $T>0$. The trajectory $t\mapsto x(t)$ being absolutely continuous, the point $x(t)$ belongs to $F$ for all $t$, so that $C_0x(t)\wedge C_1x(t)=0$, and by differentiation
$$
\dt C_0x(t)\wedge C_1x(t)=C_0A_{\lambda(t)}x(t)\wedge C_1x(t)+C_0x(t)\wedge C_1A_{\lambda(t)}x(t)=0 \quad\mbox{a.e.}
$$
that is $A_{\lambda(t)}x(t)$ is tangent to $F$ (in a weak sense because $F$ need not be regular at every point).

Let $G$ stand for the set of points $x\in \mathcal{S}^{k-1}\cap F$ that verify:
$$
\left\{
\begin{array}{ll}
x\in N \\
\mbox{ or}& \\
x\in F_0 & \mbox{and } C_0A_{\lambda(x)}x\wedge C_1x+C_0x\wedge C_1A_{\lambda(x)}x=0
\end{array}
\right.
$$
It is clear that $x(t)\in G$ for almost every $t$. Assume the set $G$ to be discrete, then the trajectory is reduced to a point, and $\lambda(t)$ is almost everywhere equal to a constant. In view of Remark (*), we have proved:

\newtheorem{instantane}{Proposition}
\begin{instantane} \label{instantane}
If the pair $(C_{\lambda},A_{\lambda})$ is observable for every $\lambda \in [0,1]$ and the set $G$ is discrete then $(\Sigma)$ is uniformly observable on $[0,T]$ for all $T>0$.

This is in particular true if $\ker C_{\lambda}=\{0\}$ for $\lambda \in[0,1]$.
\end{instantane}


\subsection{Uniform observability on $[0,+\infty)$}\label{obsinfty}

It may happen that $G$ is not discrete, though the pair $(C_{\lambda},A_{\lambda})$ is observable for every $\lambda \in [0,1]$. For instance when $\dim \K^{\bot}=1$ the condition $C_0x\wedge C_1x=0$ is void and $F$ contains an open subset of $\K$. The interior of the set $F_0$ is not empty either and the ouput of the analytic system
$$
\left\{
\begin{array}{rl}
\dot{x} & =A_{\lambda(x)} x\\
y & =C_{\lambda(x)} x
\end{array}
\right.
$$
vanishes \textbf{as long as the trajectory remains in $F_0$}. Consequently $(\Sigma)$ cannot be uniformly observable on small time intervals. However it may happen that \textbf{under the assumption that the pair $(C_{\lambda},A_{\lambda})$ is observable for every $\lambda \in [0,1]$}, no trajectory remains in $F$. We present below a proof in the case where $\dim \K$ is $1$ or $2$.

\subsubsection{$\dim \K=1$}

The sphere  $\mathcal{S}^{k-1}$ consists of two points and under the condition that the pair $(C_{\lambda},A_{\lambda})$ is observable for every $\lambda \in [0,1]$, no bad trajectory can exist (see Remark (*)).

\subsubsection{$\dim \K=2$}

We can assume without loss of generality that the rank of $C_{\lambda}$ is equal to $1$ for every $\lambda\in[0,1]$. Indeed if it vanishes for some $\lambda_0$, then the pair $(C_{\lambda_0},A_{\lambda_0})$ is not observable. If the rank of $C_{\lambda}$ is greater than $1$ for one $\lambda$ then it is greater than $1$ for all $\lambda$ except for isolated values (because the minors of the matrix $(C_0x\ C_1x)$ are polynomials of the entries), the bad trajectories are obtained for these constant inputs, and the conclusion comes from Remark (*).

In this setting we have two cases:
\begin{enumerate}
	\item $\ker C_0=\ker C_1$. Then $C_1=\alpha C_0$ for some $\alpha >0$ (if $\alpha \leq 0$ then $C_{\lambda}$ vanishes for some $\lambda\in[0,1]$).
	A bad trajectory is contained in the intersection of the one-dimensional space $\ker C_0$ with the one-dimensional sphere $\mathcal{S}^1$ and is reduced to a point.
	\item $\ker C_0\neq \ker C_1$. The set $F$ is the cone $\{(C_0x)(C_1x)\leq 0\}$. On the other hand the matrices $\exp(tA_i)$ are rotation ones ($i=0,1$). If they have the same direction of rotation, all trajectories run through the whole circle and go out of $F$. If their directions of rotation are opposite or if one is zero, then $A_{\lambda}$ vanishes for some $\lambda$. The pair $(C_{\lambda},A_{\lambda})$ is not observable for that value (see also Example \ref{Kdeux}).
\end{enumerate}

We have proved

\newtheorem{Kpetit}[instantane]{Proposition}
\begin{Kpetit} \label{Kpetit}
If $\dim \K\leq 2$ then $(\Sigma)$ is uniformly observable on $[0,+\infty)$ if and only if the pair $(C_{\lambda},A_{\lambda})$ is observable for every $\lambda \in [0,1]$.
\end{Kpetit}


\subsection{Concluding Theorem and Conjecture}\label{concluding}

We keep the notations of the previous sections. In view of Theorem \ref{CNS} and Propositions \ref{instantane} and \ref{Kpetit}, we can state:

\newtheorem{Final}[cnsHurwitz]{Theorem}
\begin{Final} \label{Final}
The switched system is GUAS as soon as the pair $(C_{\lambda},A_{\lambda})$ is observable for every $\lambda \in [0,1]$, and one of the following conditions holds:
\begin{enumerate}
	\item the set $G$ is discrete;
	\item $\dim \K\leq 2$.
\end{enumerate}
In particular the switched system is GUAS if $\ker C_{\lambda}=\{0\}$ for every $\lambda \in[0,1]$.
\end{Final}

In this theorem, only sufficient conditions are stated. However we know of no system which is not GUAS and such that the pair $(C_{\lambda},A_{\lambda})$ is observable for every $\lambda \in [0,1]$. We therefore make the following conjecture:

\noindent{\bf Conjecture}

\textit{The switched system is GUAS if and only if the pair $(C_{\lambda},A_{\lambda})$ is observable for every $\lambda \in [0,1]$.}

\vskip 0.2cm

This conjecture is equivalent to saying that the matrix $B_{\lambda}$ is Hurwitz for every $\lambda \in [0,1]$. It is possible to find systems that satisfy this condition and that are not GUAS (see \cite{Liberzon}), but we do not know such systems that satisfy the additional condition (\ref{condition}), i.e. that have a common quadratic Lyapunov function. Notice that according to our results, the dimension of the state space of such a system should be at least 4 (it should verify $\displaystyle\dim\K\geq 3$).


\subsection{Examples}\label{exemples}

\subsubsection{Hurwitz matrices} Consider the matrix
$$
\begin{pmatrix}A&-C^T\\C&D\end{pmatrix} \quad \mbox{where}\quad A=\begin{pmatrix}0&1\\-1&0\end{pmatrix},
$$
and the first line of $C$ is $\begin{pmatrix}1&0\end{pmatrix}$. According to Theorem \ref{cnsHurwitz} this matrix is Hurwitz as soon as $D$ satisfies $D^T+D<0$.

\subsubsection{Two general examples}

Let us choose a skew-symmetric $k\times k$ matrix $A$, and a $k'\times k$ matrix $C$ such that the pair $(C,A)$ is observable.
Then for any matrices $D_0$ and $D_1$ such that $D_i^T+D_i<0$ the system $\{B_0,B_1\}$ is GUAS, where:
$$
B_0=\begin{pmatrix}A&-C^T\\C&D_0\end{pmatrix} \qquad B_1=\begin{pmatrix}A&-C^T\\C&D_1\end{pmatrix}
$$
Indeed the system is in the canonical form of Lemma \ref{noyau}, and $(\Sigma)$ does not depend on $\lambda$. It is therefore uniformly observable.

In the same way, and as a direct application of Proposition \ref{instantane}, we can consider the case where the dynamics of $(\Sigma)$ is null, that is the system is defined by
$$
B_0=\begin{pmatrix}0&-C_0^T\\C_0&D_0\end{pmatrix} \qquad B_1=\begin{pmatrix}0&-C_1^T\\C_1&D_1\end{pmatrix}\quad\mbox{with}\quad    D_i^T+D_i<0.
$$
It is GUAS if and only if $C_{\lambda}$ is one-to-one for all $\lambda\in[0,1]$.


\subsubsection{An example with $\dim \K=2$}\label{Kdeux}

Consider the case where
$$
A_0=\begin{pmatrix}0&a\\-a&0\end{pmatrix}, \qquad \mbox{the first line of $C_0$ is }\ \begin{pmatrix}1&0\end{pmatrix},
$$
$$
A_1=\begin{pmatrix}0&b\\-b&0\end{pmatrix}, \qquad \mbox{the first line of $C_1$ is }\ \begin{pmatrix}0&1\end{pmatrix},
$$
the other lines of $C_0$ and $C_1$ are zero, and $D_0$, $D_1$ are $k'\times k'$ matrices, $k'\geq 1$, with $D_i^T+D_i<0$, $i=1,2$.

A straightforward computation shows that the pair $(C_{\lambda},A_{\lambda})$ is observable for every $\lambda \in [0,1]$ if and only if $a$ and $b$ are both positive or both negative: the determinant of the observability matrix is equal to\\ $\displaystyle(2\lambda^2-2\lambda+1)((1-\lambda)a+\lambda b)$.

The cone $F$ is here the set $\{(x_1,x_2)\in \R^2;\ x_1x_2\leq 0\}$, that is the union of the two orthants $\{x_1\geq 0;\ x_2\leq 0\}$ and $\{x_1\leq 0;\ x_2\geq 0\}$, and $F_0$ is equal to $F$ minus the origin. For $x=(x_1,x_2)\in F_0$ we can define
$$
\lambda(x)=\frac{-x_2}{x_1-x_2}.
$$
This system cannot be uniformly observable on small time intervals since for the feedback $x\longmapsto\lambda(x)$ the output vanishes as long as $x(t)$ belongs to $F_0$ whose interior is not empty.

However $(\Sigma)$ is uniformly observable on $[0,+\infty)$, under the condition $a$ and $b$ both positive or both negative: indeed a trajectory starting at $x\neq 0$ runs through the whole circle with radius $\norm{x}$, hence goes out of $F$.

Finally the switched system is GUAS if and only if $ab>0$.


\subsubsection{The $\dim \K=d-1$ case}\label{dmoinsun}
Let us begin by a very simple example. 
Let $A$ a be $(d-1)\times(d-1)$ skew-symmetric matrix and $C$ a $1\times(d-1)$ matrix such that 
the pair $(C,A)$ is observable, and let $d_0$ and $d_1$ be two different positive numbers. The matrices
$$
B_0=\begin{pmatrix}A&-C^T\\C&-d_0\end{pmatrix}, \qquad B_1=\begin{pmatrix}A&-C^T\\C&-d_1\end{pmatrix}
$$
define a GUAS switched system with $\dim \K=d-1$.

For a less trivial example consider the skew-symmetric $2q\times 2q$ matrix $A$ which has $q$ blocks
$$
\begin{pmatrix}0&-a_j\\a_j&0\end{pmatrix}
$$
on the diagonal and vanishes elsewhere, and
$$
C_0=\begin{pmatrix}1&0&1&0&\dots&1&0\end{pmatrix}, \qquad C_1=\begin{pmatrix}0&1&0&1&\dots&0&1\end{pmatrix}.
$$
Assume $(a_1,\dots,a_q)$ to be rationally independant. Then the orbit of $\dot{x}=Ax$ for a non zero initial state $(x^0_1,\dots,x^0_{2q})$ is dense in the torus
$$
x_{2j-1}^2+x_{2j}^2=(x^0_{2j-1})^2+(x^0_{2j})^2=T_j^2 \qquad j=1,\dots,q
$$
where at least one $T_j$ does not vanish.

Therefore this orbit meets the subset of the orthant $\{x_i\geq 0;\ i=1,\dots 2q\}$ where $x_{2j-1}>0$ and $x_{2j}>0$ for at least one $j$. But in this subset we have $(C_0x)(C_1x)>0$. This shows that every non zero orbit goes out of $F$ and that the bilinear system defined by $A_0=A_1=A$, $C_0$ and $C_1$ is uniformly observable on $[0,+\infty)$. Finally the switched system defined by the matrices
$$
B_0=\begin{pmatrix}A&-C_0^T\\C_0&-d_0\end{pmatrix}, \qquad B_1=\begin{pmatrix}A&-C_1^T\\C_1&-d_1\end{pmatrix}
$$
is GUAS for any choice of positive numbers $d_0$ and $d_1$.


\subsubsection{A singular case of the Dayawansa-Martin example}\label{Paolo}

It is a well known fact that a GUAS linear switched system has always a common strict Lyapunov function, but not always a quadratic one: in \cite{dw} Dayawansa and Martin provide an example to show that even for planar switched linear systems, GUAS does not imply the existence of a common strict quadratic Lyapunov function.

We give here an example, due to Paolo Mason, which shows that for a linear switched system, GUAS and the existence of a common non strict quadratic Lyapunov function do not either imply the existence of a common strict quadratic Lyapunov function.

Consider the $2\times 2$ switched system defined by the matrices
$$
B_0=\begin{pmatrix}-1&-1\\1&-1\end{pmatrix} \quad \mbox{ and }\quad B_1=\begin{pmatrix}-1&-3-2\sqrt{2}\\3-2\sqrt{2}&-1\end{pmatrix}.
$$  
The symmetric positive matrix
$$
P=\begin{pmatrix}1&0\\0&3+2\sqrt{2}\end{pmatrix}
$$
defines a weak quadratic Lyapunov function for this system, that is $B_i^TP+PB_i\leq 0$ for $i=0,1$.

On the other hand the switched system is GUAS: to see it, just apply Theorem 1 of \cite{bbm}. Indeed our system is in the class  satisfying the $S_4$-GUAS condition of this theorem. It remains to show that it admits no strict quadratic Lyapunov function.

We are  seeking  a positive definite symmetric
matrix P in the form
$$\left(
                         \begin{array}{cc}
                           1 & q \\
                           q & r \\
                         \end{array}
                       \right)
                        $$
such that 
\begin{equation} M_i=B_i^TP+PB_i< 0, \ \ i =1,2  \label{e1} \end{equation} 
Equation (\ref{e1}) is satisfied if the interior of the ellipses in the $(q,r)$ plan given by $\det{M_i}=0$ intersect. It is straightforward to check that those ellipses have the same major axis $q=0$ and have respectively the vertices
$$
\{(0,3 - 2 \sqrt{2}),(0, 3 + 2 \sqrt{2})\} \ \ \text{and} \ \ \{(0,3 + 2 \sqrt{2}),(0,99 + 70 \sqrt{2})\}.
$$

Consequently their interiors do not intersect.


\section{Nonlinear switched sytems}\label{nonlinear}
\subsection{The nonlinear setting}

The first part of the paper dealt with linear systems and quadratic Lyapunov functions.
The aim of this section is to extend the previous results to nonlinear switched systems admitting a weak common Lyapunov function, no longer required to be quadratic.\\

We consider the nonlinear switched system
$$
(S)\qquad \dot{X}=(1-u)f_0(X)+uf_1(X), \quad X\in\mathbb{R}^d,\quad u\in\{0,1\},
$$
where $f_0$ and $f_1$ are analytic globally asymptotically stable vector fields, and $V:\mathbb{R}^d \longrightarrow \mathbb{R}_+$ is an analytic weak common Lyapunov function, not necessarily quadratic, for the vector fields $f_0$ and $f_1$ ($V$ is a positive definite function and for $i=0,1$, the Lie derivative $\mathcal{L}_{f_i} V$ is nonpositive, so that $V$ is nonincreasing along the solutions of $(S)$). This Lyapunov function is also assumed to be radially unbounded, that is $V(X)$ goes to $+\infty$ as $\|X\|$ goes to $+\infty$, to ensure that the solutions of $(S)$ are defined on $[0,+\infty)$ for any input and initial condition.

In the same way as in the linear case we define $\mathcal{K}_i=\{\mathcal{L}_{f_i} V=0\}$, $i=0,1$, and 
the set $\mathcal{K}=\mathcal{K}_0\cap\mathcal{K}_1$. The analycity assumption is here crucial because it allows to state the analog of Lemma \ref{Inter}, that is:

\newtheorem{InterNL}[noyau]{Lemma}
\begin{InterNL} \label{InterNL}
For any input $u$ and any initial condition $X\in\mathbb{R}^d$, the $\omega$-limit set $\Omega_u(X)$ of the corresponding solution $\Phi_u(X)$ is included in $\mathcal{K}$.
\end{InterNL} 

The proof of this lemma follows the same lines as the one of Lemma~\ref{Inter} (see also \cite{JN13}). It is worthwhile to notice that is not true in general without the analycity hypothesis (a counter-example can be found in \cite{JN13}), but that the remainder of the paper makes no use of this hypothesis, so that all the forthcoming results apply to smooth systems provided that they satisfy the conclusion of Lemma~\ref{InterNL}. 

\vskip 0.2cm

Our purpose is now to define an observed system on $\K$, parametrized by $\lambda\in[0,1]$, that is to project the convexified system $\dot{X}=(1-\lambda)f_0(X)+\lambda f_1(X), \quad \lambda\in[0,1]$ onto $\K$. We do not require $\K$ to be a manifold, in particular because it is generally not at the origin (it may of course fail to be at other points). The set $\K$ being analytic is locally a finite union of analytic manifolds and we could define its tangent space at a point as the tangent space to the only such submanifold to which this point belongs. However we prefer to define the tangent space at $x\in\K$ as the linear span of the derivatives $\gamma'(0)$ of all absolutely continuous curves
 $\gamma$ contained in $\mathcal{K}$ and such that $\gamma(0)=x$. Clearly the tangent spaces at different points need not have the same dimension.

\vskip 0.2cm

The next task is to define the projection of $f_\lambda(x)=(1-\lambda)f_0+\lambda f_1$ on $T_x \mathcal{K}$. 

In the linear case $\mathcal{K}$ is a linear subspace and according to the decomposition $\mathbb{R}^d=\mathcal{K}\oplus\mathcal{K}^\perp$, the linear vector fields in $\R^d$ and the controlled and observed system on $\K$ write
$$
f_\lambda\equiv\begin{pmatrix}A_{\lambda}&-C_{\lambda}^T\\C_{\lambda}&D_{\lambda}\end{pmatrix}\quad \mbox{and}\quad (\Sigma)\ \left\{
\begin{array}{rl}
\dot{x} & =A_{\lambda} x\\
y & =C_{\lambda} x
\end{array}
\right.
$$
In other words the vector field $A_\lambda x$ is the orthogonal projection of $f_\lambda(x)$ on $\mathcal{K}$ and the component $C_\lambda x=f_\lambda(x)-A_\lambda x$ is the one on $\mathcal{K}^\perp$, the Euclidean structure being induced by the quadratic Lyapunov function.\\


	

\vskip 0.2cm

Since in the general case $V$ is nonquadratic, it does not define an Euclidean structure. However we can endow $\mathbb{R}^d$ with its canonical inner product, denoted by $\langle \cdot, \cdot\rangle$, or actually with any other inner product, it does not matter as we will see further. Consequently we can define at each point $x\in\K$ the orthogonal complement $N_x \mathcal{K}$ to $T_x\mathcal{K}$, and $g_\lambda(x)$, (respectively $h_\lambda(x)=f_\lambda(x)-g_\lambda(x)$), as the orthogonal projection of $f_\lambda(x)$ onto $T_x\mathcal{K}$, (respectively on $N_x \mathcal{K}$).

This way we get the desired observed system on $\K$
$$
(\Sigma)\ \ \  \left\{
\begin{array}{ll}
\dot{x}& =g_\lambda(x)\\
y& =h_\lambda(x)\\
\end{array}
\right.
$$ 
Of course the vector fields $g_\lambda$ are not defined in a usual way, but the solutions to $(\Sigma)$ for a given measurable input $\lambda(t)$ will merely be the absolutely continuous curves $\gamma$ in $\K$ that satisfy $\dot{\gamma(t)} =g_\lambda(\gamma(t))$ almost everywhere. On the other hand the sets $N_x \mathcal{K}$ being canonically identified with linear subspaces of $\mathbb{R}^d$ the output can be viewed as taking its values in $\mathbb{R}^d$.

To finish we need a relevant notion of observability. For linear systems, observability is equivalent to distinguishability of the origin. Of course in the non linear case distinguishing the origin of $\R^d$ does not imply observability but this weak notion will turn out to be the right one here.

\newtheorem{inftyDisting}[GUAS]{Definition}
\begin{inftyDisting}\label{inftyDisting}
The system $(\Sigma)$ is said to uniformly distinguish $0$ on $[0,+\infty)$ if for any measurable input $t\longmapsto \lambda(t)$ from $[0,+\infty)$ into $[0,1]$ for which the solution is well-defined on $[0,+\infty)$ , the output distinguishes the origin. 
\\A  solution $x(t)\neq0$ of $(\Sigma)$ defined on $[0,+\infty)$ for which the output vanishes is said to be a bad trajectory. 	
\end{inftyDisting}

\vskip 0.2cm

 \noindent{\bf Remark}
A solution which is not defined for all positive times reaches the boundary of $\mathcal{K}$ in a transverse way, actually leaves $\K$ (recall that all solutions of $(S)$ are defined in positive time). It is consequently not a limit trajectory of the initial switched system, and is not a ''bad trajectory'' (in the same sense as in the linear case).

\vskip 0.5cm

We are now in a position to state the main result of this section.

\newtheorem{GUASnon}[cnsHurwitz]{Theorem}
\begin{GUASnon}\label{GUASnon}
The switched system $(S)$ is GUAS if and only if $(\Sigma)$ uniformly distinguishes $0$ on $[0,+\infty)$.
\end{GUASnon}

\demo
The proof of the theorem follows the same steps as the one of Theorem~\ref{CNS}. The tools are a bit different and come in part from \cite{JN13}. \\
 Let us assume that the switched system is not GUAS. Then there exist a measurable input $u:[0,+\infty) \longrightarrow \{0,1\}$ and an initial state $X\in\mathbb{R}^d$ for which the solution $\Phi_u(t,X)$ does not converge to the origin.\\
Let $l\neq0$ be an $\omega$-limit point of $\Phi_u(t,X)$, and $(t_j)_{j\geq0}\subset \mathbb{R}_+$ be an increasing sequence such that
 	$$\lim_{j\to +\infty} \Phi_u(t_j,X)=l.$$
Consider the sequence $(\phi_j)_{j\geq0}$ defined by
	$$\phi_j(t)=\Phi_u(t_j+t,X), \forall t\geq0, j\geq0.$$ 
Following the same lines as in the linear case (see also \cite{JN13}) we can prove:
\newtheorem{Convergence}[noyau]{Lemma}
\begin{Convergence} \label{Converegnce}
Up to a subsequence, $(\phi_j)_{j\geq0}$ converges uniformly on each compact to an absolutely continuous function $\phi:[0,+\infty) \longrightarrow \Omega_u(X)$.\\
Moreover $\phi(0)=l$ and there exists a measurable function $\lambda:[0,+\infty) \longrightarrow [0,1]$ such that
	$$\dot{\phi}=(1-\lambda)f_0(\phi)+\lambda f_1(\phi) \mbox{ a.e.}$$
In other words, $\phi$ is a solution of the convexified switched system i.e. $\dot{\phi}=f_\lambda(\phi)$ a.e.	
\end{Convergence} 

Thanks to Lemma~\ref{InterNL}, for all $t\in[0,+\infty)$, $\phi(t)\in\mathcal{K}$. It follows that, for almost all $t\in[0,+\infty)$, $\dot{\phi}(t)\in T_{\phi(t)} \mathcal{K}$
and 
$\left\{ 
\begin{array}{rcl}
g_\lambda(\phi)&=&f_\lambda(\phi) \\
h_\lambda(\phi)&=&0
\end{array}
\right.$ a.e.\\
Therefore $\phi\neq0$ is a solution of the controlled and observed system $(\Sigma)$ for which the output vanishes.\\
Consequently $(\Sigma)$ does not uniformly distinguish $0$ on $[0,+\infty)$.\\

Conversely, let us assume that the controlled and observed system $(\Sigma)$ does not uniformly distinguish $0$.\\
There exists a measurable input $\lambda: [0,+\infty) \longrightarrow [0,1]$ such that the trajectory $\phi$ for some initial condition $\phi(0)\neq0$ satisfies $h_\lambda(\phi)=0$ almost everywhere.
Then $f_\lambda(\phi)=g_\lambda(\phi)$ almost everywhere and, since $\phi(t)\in\K=\K_0\cap\K_1$, 
$$
\dt V(\phi(t))=(1-\lambda(t))L_{f_0}V(\phi(t))+\lambda(t)L_{f_1}V(\phi(t))=0 \quad a.e.
$$
This shows that $\phi$ is contained in a level set of $V$, and is a trajectory of the convexified switched system that does not converge to the origin.\\
In conclusion the convexified switched system is not GUAS and neither is the switched system (S).

\hfill  $\Box$\\

\vskip 0.2cm
\noindent{\bf Remark}.
This proof shows clearly that the results do not depend on the choice of the inner product (or even on a Riemannian structure) on $\R^d$. Indeed, the output $h_\lambda$ is solely required to vanish if and only if the vector field $f_\lambda$ is tangent to the analytic set $\mathcal{K}$, and this do not depend on the chosen inner product.






\subsection{Applications}

The nonlinear counterpart of Proposition~\ref{Kpetit} is no longer true. However for the same reasons as for linear systems, if the pair $(g_\lambda,h_\lambda)$ distinguishes $0$ for each $\lambda\in[0,1]$ then no bad trajectory can be constant. Hence, Proposition~\ref{Kpetit} remains true if $\dim \mathcal{K}=0$, that is if $\K$ is discrete. \\

Moreover Conditions guaranteeing that $(\Sigma)$ uniformly distinguishes $0$ on $[0,+\infty)$ can be derived from Theorem~\ref{GUASnon} in the case when $\mathcal{K}$ is a one-dimensional submanifold, at least outside of the origin. 
Namely, a necessary and sufficient condition for $(\Sigma)$ to uniformly distinguish $0$ on $[0,+\infty)$ is that no bad periodic trajectory exists. \\

Let us assume that $\mathcal{K}\setminus \{0\}$ is a one-dimensional manifold and let $x(t)$ be a bad trajectory. It lies in a level set $\{V=r\}$ of the Lyapunov function. Since the set $\mathcal{K}\cap\{V=r\}$ is compact, it can be covered by a finite number of coordinate charts $(U_i,\psi_i)_{1\leq i\leq n}$ such that $\psi_i(\mathcal{K}\cap U_i)$ is an open interval for $i=1,\ldots,n$.\\
All we need to conclude lies in the forthcoming remarks:\\

 \noindent{\bf Remarks}
\begin{enumerate}
	\item Since the pair $(g_\lambda,h_\lambda)$ distinguishes $0$ for each $\lambda\in[0,1]$ no bad trajectory can be constant. Therefore the trajectory $x(t)$ cannot converge to a point.
	\item It follows from the first item that for each $\lambda\in[0,1]$, the vector field $g_\lambda$ does not vanish on the set $N=\{x\in\mathcal{K}: \quad h_0(x)=h_1(x)=0\}$. This is equivalent to saying that the vector fields $g_0$ 		and $g_1$ point in the same direction on $N$. 
	\item For $x\in\mathcal{K}$ and $\lambda\in[0,1]$, $h_\lambda(x)=0$ implies that either $\lambda=\lambda(x)$ is unique or $x\in N$. \\
		Up to a reparametrization of time 
		$\dot{x}=\left\{ 
		\begin{array}{l}
			g_0(x) \mbox{\quad on } N \\
			g_{\lambda(x)}(x) \mbox{\quad elsewhere}
		\end{array}
		\right.$.\\
		The solution $x(t)$ goes through $\mathcal{K}$ in the same direction. In other words no turn around is allowed.
		
\end{enumerate}
 
\noindent As a consequence the solution $x(t)$ leaves the coordinate charts in finite time so it has to come back to one it already has visited.
This ensures that there exists a positive time $t$ such that $x(t)=x(0)$ and proves the existence of a bad periodic trajectory.\\
This reasoning still works when $\dim \mathcal{K}=2$ and for each $r>0$ the connected components of $\mathcal{K}\cap\{V=r\}$ 
are submanifolds whose dimension is lower than or equal to $1$. \\

\newtheorem{KpetitNL}[instantane]{Proposition}
\begin{KpetitNL}\label{KpetitNL}
Assume that the pair $(g_\lambda,h_\lambda)$ distinguishes $0$ for each $\lambda\in[0,1]$, and that one of the following conditions hold:
\begin{enumerate}
\item $\dim \K=0$,
\item  $\mathcal{K}\setminus \{0\}$ is a one-dimensional manifold that contains no periodic trajectory of the convexified system,
\item $\dim \mathcal{K}=2$ and for each $r>0$ the connected components of $\mathcal{K}\cap\{V=r\}$ 
		are submanifolds whose dimension is $0$ or $1$, and contain no periodic trajectory of the convexified system in that last case,
\end{enumerate}
then $(\Sigma)$ uniformly distinguishes $0$ on $[0,+\infty)$, and the switched system is GUAS.
\end{KpetitNL}

As $\K$ is an analytic set, there exists a neighborhood of the origin that encounters only one conected component of $\K$, the one that contains the origin itself. Consequently we can also state the following local result.
		   					   
\newtheorem{LocalNL}[instantane]{Proposition}
\begin{LocalNL}\label{LocalNL}
Assume that the pair $(g_\lambda,h_\lambda)$ distinguishes $0$ for each $\lambda\in[0,1]$, and that there exists a neighborhood $W$ of the origin such that 
$W\bigcap\mathcal{K}\setminus \{0\}$ is a one-dimensional manifold.

Then the switched system is locally asymptotically stable.
\end{LocalNL}
\demo
We can assume first that $W\cap\K$ is connected, and second that $W$ to is of the form $\{V<r\}$ for some $r>0$, so that no trajectory starting in $W$ can leave $W$. Due to the connectedness of $W\cap\K$ no trajectory starting in $W$ can be periodic, so that we can apply the previous reasonning.

\hfill $\Box$


\subsection{Example}
 
We provide now an example where $\dim \mathcal{K}=2$ that illustrates Proposition \ref{KpetitNL}.
Consider the two analytic vector fields on $\mathbb{R}^3$ defined by

	$$f_0(x)=\begin{pmatrix} -x_2-x_3    \\
					         x_1 \\
					         x_1-x_3\end{pmatrix}$$
and
	$$f_1(x)=\begin{pmatrix} -x_2-\varphi(x_1,x_2)x_3    \\
					         x_1 \\
					         \varphi(x_1,x_2)x_1-x_3 \end{pmatrix}$$
where $\varphi:\mathbb{R}^2 \longrightarrow \mathbb{R}$ is an analytic function which is non constant on spheres.\\
Consider the weak Lyapunov function $V(x)=x_1^2+x_2^2+x_3^2$. A straight computation gives $L_{f_0}V=L_{f_1}V=-2x^3$, so that $\mathcal{K}=\{x_3=0\}$.\\
For each $x=(x_1,x_2)\in\mathcal{K}\equiv\mathbb{R}^2$	, 	

	$$f_0(x)=\begin{pmatrix} -x_2    \\
					         x_1 \\
					         x_1\end{pmatrix}
					         				         \quad\mbox{and}\quad
	f_1(x)=\begin{pmatrix} -x_2  \\
					         x_1 \\
					         \varphi(x)x_1 \end{pmatrix}
					         $$		         					         
and the observed system on $\K$ is: 					         

$$
\left\{
\begin{array}{ll}
g_\lambda(x)& =\begin{pmatrix} -x_2  \\
					         x_1  \end{pmatrix}\\
h_\lambda(x)&=(1-\lambda+\lambda\varphi(x))x_1\\
\end{array}.
\right.
$$
Its trajectories do not depend on the input and are circles around the origin. Since $\varphi$ is non constant on spheres, $(\Sigma)$
distinguishes $0$ on $[0,+\infty)$  for each constant $\lambda\in[0,1]$ (it actually distinguishes $0$ on $[0,2\pi)$). 

Let us first assume that $\varphi$ takes its values in $\mathbb{R}_{-}^*$. The feedback $\lambda=\dfrac{1}{1-\varphi}\in(0,1)$ is well defined in that case, and as $1-\lambda+\lambda\varphi=0$, we get  $h_{\lambda(x)}(x)=0$ for each $x\in\mathbb{R}^2$. Consequently, the input $\lambda$ does not distinguish $0$
on $[0,+\infty)$ and according to Theorem~\ref{GUASnon} the switched system is not GUAS.

On the other hand, if the function $\varphi$ is allowed to take positive values on each sphere, for instance $\varphi(x)=\frac{1}{2}\sin x_1$, then no bad trajectory can exist. 
Indeed, all trajectories are periodic and must go through a point $x$ for which $\varphi(x)$ is positive. Hence, $h_0$ and $h_1$ do not vanish and have the same sign
in a neighborhood of $x$. Consequently the output cannot vanish in this neighborhood. According to Theorem~\ref{GUASnon} the switched system is GUAS.


\section{Conclusion}\label{Baratin}

In this work we have presented a new method to obtain sufficient conditions of global asymptotic stability of pairs of Hurwitz matrices that share a weak quadratic Lyapunov function. We have also shown that these linear results partially extend to pairs of nonlinear analytic vector fields that share a non quadratic weak Lyapunov function.

A challenging question is now the observability property on $[0,+\infty)$ of the bilinear systems we have encountered, in the case when the dimension of $\K$ is greater than $2$.

On the other hand our results can be easily extended to finite families of matrices or vector fields, but the result is rather tedious to write down. In the linear case let $\{B_1,\dots,B_p\}$, with $p\geq 3$, be a finite family of Hurwitz matrices that share a weak quadratic Lyapunov function, that can be assumed to be the identity. We claim that the switched system they define is GUAS if and only if:
\begin{enumerate}
	\item[($P_2$)] For each pair $i\neq j\in \{1,\dots,p\}$ the observed system on $\K_i\bigcap\K_j$ is observable on $[0,+\infty)$;
	\item[($P_3$)] Property $P_2$ holds and for each 3-uple $i,j,k\in \{1,\dots,p\}$ the observed system on $\K_i\bigcap\K_j\bigcap\K_k$ is observable on $[0,+\infty)$;
	\item[($P_k$)] and so on, up to
	\item[($P_p$)] Properties $P_2$ to $P_{p-1}$ hold and the observed system on $\bigcap_{i=1}^p\K_i$ is observable on $[0,+\infty)$.
\end{enumerate}

The key point of this statement lies in the proof of Lemma \ref{Inter}.


\vskip 0.3cm

\noindent{\bf Acknowledgments}. The authors wish to express their thanks to Paolo Mason for the example of Section \ref{Paolo}.


\end{document}